\documentclass[twocolumn]{amsart}
\pdfoutput=1
\usepackage{latexsym, amsmath, amssymb, amsthm}
\usepackage{rsfso}
\usepackage{mathtools}
\usepackage{paralist}
\usepackage{amsrefs}
\usepackage[T1]{fontenc}
\usepackage{mathpazo}
\usepackage{microtype}
\usepackage{mathrsfs}
\usepackage{csquotes}
\usepackage{tikz}
\usepackage[colorlinks=true, linkcolor=blue, urlcolor=blue, citecolor=blue]%
    {hyperref}
\usepackage{cleveref}
\usepackage[centering, includeheadfoot, hmargin=0.75in, tmargin=0.6in, 
  bmargin=0.75in, headheight=6pt]{geometry}
\newtheorem{lemma}{Lemma}[section]
\newtheorem{theorem}[lemma]{Theorem}

\theoremstyle{definition}

\newtheorem{question}[lemma]{Question}
\newtheorem{examplex}[lemma]{Example}
\newenvironment{example}
  {\pushQED{\qed}\examplex}
  {\popQED\endexamplex}

\newcommand{\relphantom}[1]{\mathrel{\phantom{#1}}}

\newcommand{\CC}{\ensuremath{\mathbb{C}}} 
\newcommand{\PP}{\ensuremath{\mathbb{P}}} 

\newcommand{\RR}{\ensuremath{\mathbb{R}}} 
 
\DeclareMathOperator{\codim}{codim}
\DeclareMathOperator{\conv}{conv}

\DeclareMathOperator{\genus}{p_{\text{a}}}
\DeclareMathOperator{\hilbreg}{m}

\DeclareMathOperator{\py}{py}

\DeclareMathOperator{\qp}{qp}

\begin{document}

\title[Sums of squares: a real projective story]%
{Sums of squares: \\[-2pt] a real projective story}

\author[G.~Blekherman]{Grigoriy Blekherman}
\address{Grigoriy Blekherman: School of Mathematics, Georgia Institute of Technology, 686 Cherry
  Street, Atlanta, Georgia, 30332, United States of America;
  {\normalfont \texttt{greg@math.gatech.edu}}}

\author[R.~Sinn]{Rainer Sinn}
\address{Rainer Sinn: Mathematisches Institut, Augustusplatz 10, 
  Universit\"at Leipzig, 04109 Leipzig, Germany; {\normalfont
    \texttt{rainer.sinn@uni-leipzig.de}}}

\author[G.G.~Smith]{Gregory G.{} Smith}
\address{Gregory G.{} Smith: Department of Mathematics and Statistics, Queen's
  University, Kingston, Ontario, K7L 3N6, Canada; {\normalfont
    \texttt{ggsmith@mast.queensu.ca}}}

\author[M.~Velasco]{Mauricio Velasco} 
\address{Mauricio Velasco: Departamento de Matem\'aticas, Universidad de los
  Andes, Carrera 1 No. 18a 10, Edificio~H, Primer Piso, 111711 Bogot\'a,
  Colombia; {\normalfont \texttt{mvelasco@uniandes.edu.co}}}


\maketitle

\begin{figure}[ht]
  \centering
    \includegraphics[width = \the\columnwidth]{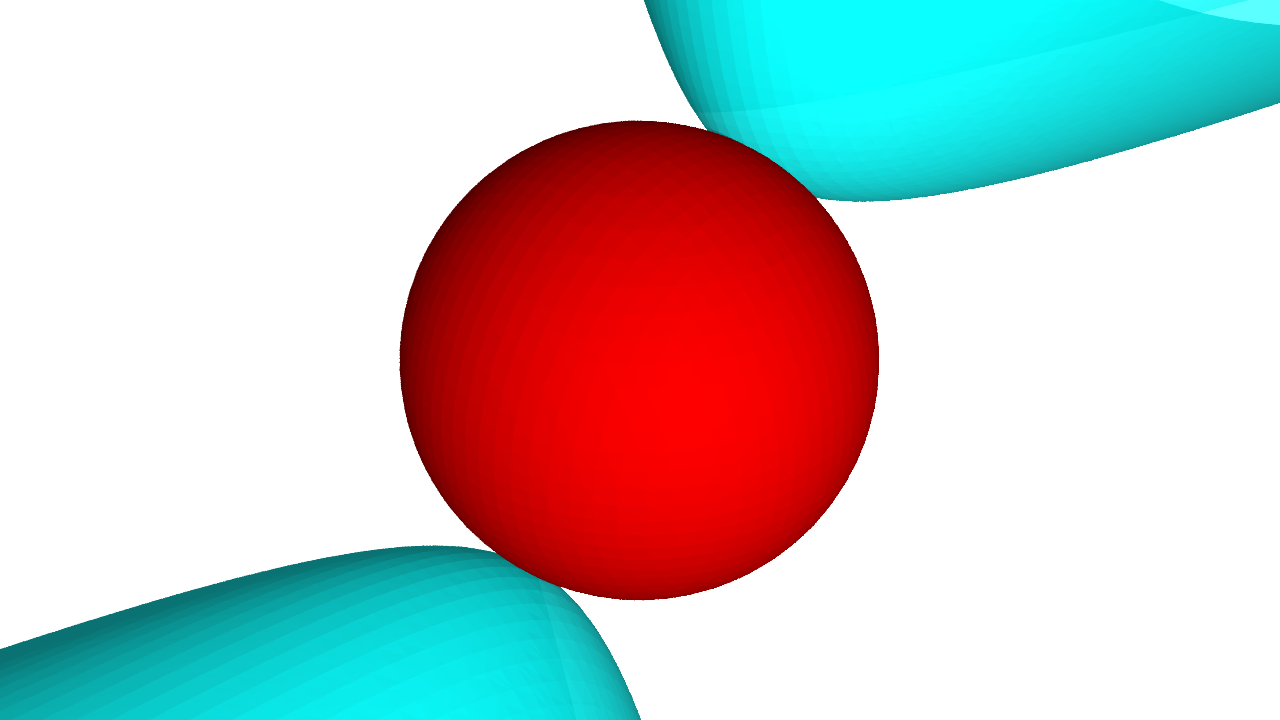}
    \caption{The unit sphere in red and the zero-locus of $2(x^2-yz)+1$ in
      teal}
  \label{f:sphere}
\end{figure}

\noindent
How do we find the minimum of the function $2(x^2 - yz)$ on the unit
sphere? An algebraic way to see that the minimum equals $-1$ involves the
decomposition
\[
  2(x^2 - yz) + 1 - ( 1-x^2-y^2-z^2) = 3 x^2 + (y-z)^2 \, .
\]
The unit sphere in $\RR^3$ is the zero-locus of the polynomial
$1-x^2-y^2-z^2$.  This decomposition establishes that the function
$2(x^2-yz)+1$ is a sum of squares modulo the defining equation of the unit
sphere, so this function is nonnegative on the sphere.  Thus, we have a
certificate that the minimum of $2(x^2-yz)$ on the unit sphere is at least
$-1$.  This lower bound is optimal because the sum of squares $3x^2+(y-z)^2$
vanishes at the points $\pm \tfrac{1}{\sqrt{2}}(0, 1, 1)$ on the unit sphere.
Figure~\ref{f:sphere} illustrates that the unit sphere is tangent to the
zero-locus of the function $2(x^2-yz)+1$ at these two points.

This example indicates that nonnegativity is intimately related to polynomial
optimization and sum-of-squares representations give rise to lower bounds on
minima.  Does this approach lead to sharp bounds for all quadratic
polynomials?  Does it apply to higher-degree polynomials on the unit sphere?
What happens when we replace the unit sphere with another algebraic variety?
This article examines these basic questions.  In the process, we also uncover
fascinating connections between real and complex algebraic geometry.
Remarkably, convex geometry provides the bridge between these two worlds.

\begin{figure}[ht]
  \centering
    \includegraphics[width = \the\columnwidth]{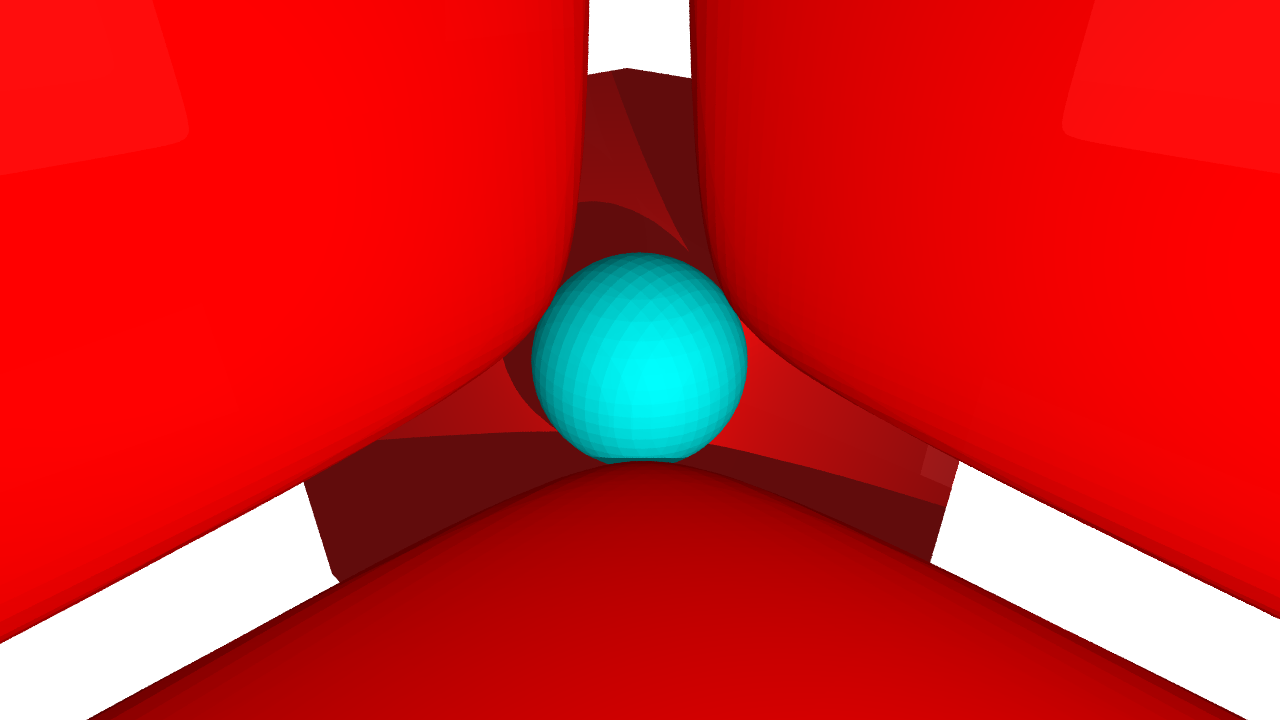}
    \caption{The real zero-loci of $xyz-1$ in red and $x^2+y^2+z^2-3$ in
      teal}
  \label{f:motzkin}
\end{figure}

As a second motivating example, let $X$ be the cubic surface defined as the
real zero-locus of the polynomial $q \coloneqq xyz - w^3$ and consider
$f \coloneqq x^2 + y^2 + z^2 - 3w^2$.  This quadratic polynomial is not a sum
of squares modulo the defining ideal of $X$ because no nonzero polynomial of
degree less than $3$ vanishes on the algebraic variety $X$ and $f$ is not
globally nonnegative.  Nevertheless, we claim that $f$ is nonnegative on the
algebraic variety $X$.  We prove this in \Cref{sec:psdsos} by recognizing $f$
as the Motzkin polynomial in disguise. More directly, we certify the
nonnegativity of $f$ on $X$ with the sum-of-squares decomposition
\[
  g \, f - 12w \, q = h_1^2 + h_2^2 + \dotsb + h_6^2
\]
where $g \coloneqq x^2 + y^2 + z^2 + (\!\sqrt{3} \, w)^2$ is a sum of squares,
and the squares on the right are given by $h_1 \coloneqq x^2-w^2$,
$h_2 \coloneqq \!\sqrt{2} (xy-zw)$, $h_3 \coloneqq \!\sqrt{2} (xz-yw)$,
$h_4 \coloneqq y^2-w^2$, $h_5 \coloneqq \!\sqrt{2} (yz-xw)$, and
$h_6 \coloneqq z^2-w^2$.  The important new ingredient is the multiplier $g$.
As $g$ is nonnegative and not divisible by the irreducible polynomial $q$,
this decomposition confirms that $f \geqslant 0$ whenever $q = 0$.
Equivalently, the function $f$ is the sum $\sum_{i} h_i^2/g$ of nonnegative
rational functions on the variety $X$.  Restricting to $w = 1$,
Figure~\ref{f:motzkin} depicts the variety $X$ as tangent to the sphere
defined by $x^2+y^2+z^2=3$ at the points $(1,1,1)$, $(1,-1,-1)$, $(-1,1,-1)$,
and $(-1,-1,1)$.  Thus, the restriction of the function $f$ to
$X \cap \{ w \neq 0 \}$ has a minimum value equal to zero.

To understand the relationship between nonnegativity and sums of squares, it
suffices to consider homogeneous polynomials.  Indeed, homogenizing any
sum-of-squares decomposition by introducing a new variable produces a
homogeneous sum-of-squares decomposition.  It follows that a polynomial is a
sum of squares modulo the ideal of an affine variety
$X \subseteq \mathbb{A}^n$ if and only if an appropriate homogenization
(having sufficiently high degree) is a sum of squares modulo the homogeneous
ideal of the projective completion $\overline{X} \subseteq \PP^n$.

Working with projective varieties has conceptual and technical advantages, so
we concentrate on nonnegative functions and sums of squares on real
subvarieties of projective space.  We examine three problems:
\begin{compactenum}[\upshape (1)]
\item Identify all of the real projective varieties on which every nonnegative
  quadratic function is a polynomial sum of squares.
\item Knowing that a nonnegative function is a polynomial sum of squares on a
  projective variety, control the number of squares needed in such a
  decomposition.
\item Bound the degree of a sum-of-squares multiplier such that its product
  with a nonnegative function decomposes into a polynomial sum of squares.
\end{compactenum}
Each of the subsequent sections is devoted to one of these problems.  Sections
begin with some historial context before describing more recent results.

Sums of squares of polynomials are indispensable in the theory of nonnegative
polynomials.  Beyond forming a self-evident subset of nonnegative polynomials,
there are practical algorithms for deciding whether a given polynomial is a
sum of squares.  The existence of a polynomial sum-of-squares decomposition is
equivalent to the feasibility of a semidefinite programming problem: a convex
optimization problem with a linear objective function on the intersection of
the cone of positive semidefinite matrices with an affine linear space.  There
are robust and efficient software packages, typically based on interior point
methods, for solving these semidefinite programming problems.  More
theoretically, polynomial sums of squares generate an infinite hierarchy of
approximations for the set of nonnegative polynomials.  The resulting
sum-of-squares relaxations (also known as Lasserre relaxations) play a
prominent role in engineering applications and computer science; see Chapter~7
in \cite{SIAMbook}.

\section{Geometry of sums of squares}
\label{sec:psdsos}

\noindent
When is every nonnegative polynomial a sum of squares?  For which degrees and
in how many variables is there an equivalence? David Hilbert gave a
eulogy~\cite{eulogy} in 1910 for Hermann Minkowski (1864--1909) to the Royal
Society of Sciences in G\"ottingen. Hilbert explained that the history of this
question began in 1885 with Minkowski's thesis defense. Hilbert participated
in this defense as an examiner. Two years prior to the defense, Minkowski made
a name for himself by solving a problem posed by Eisenstein in a competition
for the Academy of Sciences (Paris) about the number of representations of a
positive integer as a sum of five squares.  His doctoral thesis was a
continuation of this prize winning work.  In his defense, Minkowski argued
that there are polynomials with real coefficients that are nonnegative
functions on Euclidean space but cannot be written as sums of squares of
polynomials.

Minkowski never published these insights, but he did inspire Hilbert's seminal
work in the area. Hilbert~\cite{Hilbert1888} demonstrated, in 1888, that
Minkowski's claim is correct by completely characterizing when every
nonnegative polynomial is a polynomial sum of squares.
\begin{theorem}[Hilbert]
  \label{thm:hilb}  
  Every nonnegative homogenous polynomial on $\RR^n$ having degree $2d$ is a
  sum of squares of polynomials if and only if
  \begin{compactitem}[$\bullet$]
  \item $2d = 2$, \hfill \emph{(quadratic forms)}    
  \item $n = 2$, or \hfill \emph{(binary forms)}
  \item $n = 3$ and $2d = 4$. \hfill \emph{(ternary quartics)}
  \end{compactitem}
\end{theorem}

\noindent
In a second paper, Hilbert~\cite{Hilbert1893} also characterized nonnegative
polynomials in two variables (or equivalently homogeneous polynomials in three
variables) as sums of squares of rational functions.  These results
undoubtably prompted Hilbert to formulate the 17th problem on his celebrated
list of 23 problems.

For the first two cases in Theorem~\ref{thm:hilb}, establishing that all
nonnegative polynomials are sums of squares of polynomials is relatively
straightforward.  The third case is noticeably more challenging.  Arguably,
the hardest part is proving that there exist nonnegative polynomials that are
not polynomial sums of squares in all of the remaining cases. Hilbert
accomplished this task with an impressive nonconstructive argument.

Curiously, Theodore Motzkin published in 1967, nearly 80 years after Hilbert's
paper, the first explicit example of a nonnegative polynomial that cannot be
expressed as a polynomial sum of squares. The Motzkin polynomial is the
ternary sextic $x_0^6 + x_1^4 x_2^2 + x_1^2 x_2^4 - 3 x_0^2 x_1^2 x_2^2$.  By
taking suitable means of $x_0^2$, $x_1^2$, and $x_2^2$, the nonnegativity of
this homogeneous polynomial on $\RR^3$ follows immediately from the
arithmetic-geometric mean inequality.

\subsection*{Real projective varieties}

As alluded in the preceding section, it is advantageous to work with
homogenous polynomials and real projective varieties.  Within this broader
framework, we not only generalize Theorem~\ref{thm:hilb} but discover a
geometric characterization of the equality between nonnegative polynomials and
sums of squares.

We consider the $n$-dimensional projective space $\PP^n$ to be the set of
lines in $\CC^{n+1}$ passing through the origin.  The point
$p \coloneqq [p_0 \mathbin{:} p_1 \mathbin{:} \dotsb \mathbin{:}
\mathbin{p}_n] \in \PP^n$ is the equivalence class
\[
  \bigl\{ \!(\lambda p_0, \lambda p_1, \dotsc, \lambda p_n) \!\in \CC^{n+1}
  \setminus \{ \mathbf{0} \} \mid \text{for some $0 \neq \lambda \!\in \CC$}
  \bigr\} .
\]
The point $p$ is real if it has a representative
$[p_0 \mathbin{:} p_1 \mathbin{:} \dotsb \mathbin{:} \mathbin{p}_n]$ such that
$p_i \in \RR$ for all $0 \leqslant i \leqslant n$.

A real subvariety $X \subseteq \PP^n$ is the set of common zeroes of some
homogeneous polynomials with real coefficients in the variables
$x_0, x_1, \dotsc, x_n$.  Whenever a homogenous polynomial $f$ vanishes at
$p \coloneqq (p_0, p_1, \dotsc, p_n) \in \CC^{n+1}$, it also vanishes at all
scalar multiples because homogeneity implies that
$f(\lambda p) = \lambda^{\deg(f)} f(p)$ for all nonzero scalars $\lambda$.  In
other words, the homogenous polynomial $f$ vanishes at the point
$[p_0 \mathbin{:} p_1 \mathbin{:} \dotsb \mathbin{:} \mathbin{p}_n] \in \PP^n$
and the set of common zeroes of a collection of homogenous polynomials is a
subset of $\PP^n$.  For example, the twisted cubic curve
\[
  C \coloneqq \bigl\{ [t_0^3 \mathbin{:} t_0^2 t_1^{} \mathbin{:} t_0^{} t_1^2
  \mathbin{:} t_1^3 ] \in \PP^3 \mid [t_0^{} \mathbin{:} t_1^{} ] \in \PP^1
  \bigr\}
\]
is the set of common zeroes of the quadratic polynomials
$x_0^{} x_2^{} - x_1^2$, $x_0^{} x_3^{} - x_1^{} x_2^{}$, and
$x_1^{} x_3^{} - x_2^2$.

Throughout, we assume that the subvariety $X \subseteq \PP^n$ is nondegenerate
and totally-real.  The first assumption means that the subvariety $X$ does not
lie in a hyperplane.  This mild hypothesis guarantees that the ambient
projective space $\PP^n$ is not unnecessarily large.  The second assumption is
more significant---it ensures that $X$ has enough real points.  By definition,
the subvariety $X$ is totally-real if the set of real points in $X$, regarded
as a subset of the complex points in $X$, is dense in the Zariski topology.
Equivalently, every irreducible component of the algebraic variety $X$ has a
nonsingular real point.  The twisted cubic curve in $\PP^3$ is totally-real
whereas the zero set of the polynomial $x_0^2 + x_1^2 + \dotsb + x_n^2$ is not
because it does not contain a real point in $\PP^n$.

In this geometric context, polynomials are replaced by elements in another
ring.  The homogeneous coordinate ring $R$ of the subvariety
$X \subseteq \PP^n$ is the quotient of the polynomial ring
$\RR[x_0, x_1, \dotsc, x_n]$ by the ideal generated by all homogenous
polynomials that vanish on $X$.  For any nonnegative integer $j$, let $R_j$
denote the real vector space spanned by the homogeneous elements in $R$ of
degree $j$.  Generalizing the concept of a polynomial sum of squares is not
difficult.  A homogeneous element $f \in R_{2j}$ is a sum of squares on $X$ if
there exists homogeneous elements $h_1, h_2, \dotsc, h_r \in R_{j}$ such that
$f = h_1^2 + h_2^2 + \dotsb + h_r^2$.

In comparison, defining nonnegativity requires more care. A homogeneous
element $f \in R_{2j}$ of even degree is nonnegative on $X$ if its evaluation
at each real point in $X$ is greater than or equal to $0$.  Since elements in
the ring $R$ and points in the space $\PP^n$ may both be thought of as
equivalence classes, the evaluation process involves choosing a polynomial in
$\RR[x_0,x_1, \dotsc, x_n]$ to represent the ring element and a point in
$\RR^{n+1}$ to represent the real point in $X \subset \PP^n$.  Although a
homogenous element $f$ in $R_{2j}$ does not determine a function from $X$ to
$\RR$, evaluating its polynomial representative at a representative point in
$\RR^{n+1}$ does have a well-defined sign because the degree of polynomial is
even.  We recover our original notion of nonnegativity for polynomials when
$X = \PP^n$

\begin{example}
  Since we have
  \[
    x_0^2 - x_0^{} x_3^{} + x_1^{} x_2^{} + x_2^{} x_3^{}
    = x_0^2 + x_2^{} x_3^{} - (x_0^{} x_3^{} - x_1^{} x_3^{}) \, , 
  \]
  both $x_0^2 - x_0^{} x_3^{} + x_1^{} x_2^{} + x_2^{} x_3^{}$ and
  $x_0^2 + x_2^{} x_3^{}$ represent the same element in the homogeneous
  coordinate ring of the twisted cubic curve $C$.  Choosing
  $(1,0,0,0) \in \RR^{4}$ or $(-2,0,0,0) \in \RR^{4}$ as representatives for
  the real point $[1 \mathbin{:} 0 \mathbin{:} 0 \mathbin{:} 0] \in C$, either
  of the evaluations
  \[
    (1)^2 - (1)(0) + (0)(0) + (0)(0) = 1 > 0
  \]
  or $(-2)^2 + (0)(0) = 4 > 0$ show that this element in the homogeneous
  coordinate ring is positive at this point.
\end{example}

\subsection*{Sums of squares on real projective varieties}
The real vector space $R_2$ of quadrics on the real subvariety
$X \subseteq \PP^n$ contains two fundamental subsets.
\begin{compactitem}[$\bullet$]
\item The set $P_X \subseteq R_2$ consists of those elements whose evaluations
  at every real point in $X$ is nonnegative;
  \[
    P_X \coloneqq \{ f \in R_2 \mid \text{ $f(x) \geqslant 0$ for all real
      points $x \in X$} \} \, .
  \]
\item The set $\Sigma_X \subseteq R_2$ consists of sums of squares;
  \[
    \Sigma_X \coloneqq \{ \textstyle\sum_i h_i^2 \mid \text{$h_i \in R_1$ for
      all $i$} \} \, .
  \]
\end{compactitem}
As the square of any real number is nonnegative, we have the inclusion
$\Sigma_X \subseteq P_X$. Following Hilbert~\cite{Hilbert1888}, we ask the
following question.

\begin{question}
  \label{q:general}
  For which subvarieties $X \subseteq \PP^n$ is every nonnegative quadratic
  element in $R$ a sum of squares?  Equivalently, when does the equality
  $\Sigma_X = P_X$ hold?
\end{question}

At first glance, focusing on just the quadratic elements seems to be a limited
generalization of Hilbert's work.  However, the elbow room gained by
considering all real projective subvarieties alleviates this concern.  Suppose
that we are interested in the homogeneous elements of degree $2j$ on an
subvariety $X \subseteq \PP^n$.  Set $m \coloneqq \smash{\binom{n+j}{j}} - 1$
and let $\nu_j \colon \PP^n \to \PP^m$ be the $j$th Veronese map that sends
the point $[p_0 \mathbin{:} p_1 \mathbin{:} \dotsb \mathbin{:} p_n] \in \PP^n$
to the point in $\PP^m$ whose coordinates are all possible monomials of degree
$j$ evaluated at $(p_0, p_1, \dotsc, p_n)$.  By re-embedding the subvariety
$X \subseteq \PP^n$, it is enough to understand the quadratic elements on the
image $\nu_j(X) \subset \PP^m$.

\begin{example}
  \label{e:twi}
  On the twisted cubic curve $C$, the set $P_C$ may be identified with the
  homogeneous polynomials in $\RR[t_0, t_1]$ having degree $6$ that are
  nonnegative on $\PP^1$. Likewise, the set $\Sigma_X$ may be identified with
  sums of squares of homogeneous polynomials in $\RR[t_0, t_1]$ having degree
  $3$. By Theorem~\ref{thm:hilb}, these sets coincide, so $P_C = \Sigma_C$.
\end{example}

\begin{example}
  A variant of the Veronese map shows that the Motzkin polynomial
  $x_0^6 + x_1^4 x_2^2 + x_1^2 x_2^4 - 3 x_0^2 x_1^2 x_2^2$ is not a sum of
  squares of cubic polynomials.  If it were, then an easy (if somewhat
  tedious) case study shows that the cubic polynomials would only involve the
  monomials $x_0^3, x_1^2 x_2^{}, x_1^{} x_2^2, x_0^{} x_1^{} x_2^{}$.
  Alternatively, one may use Newton polytopes (the convex hull of the exponent
  vectors) to identify these monomials.  The Newton polytope of a sum of
  squares contains the Newton polytopes of every summand and the Newton
  polytope of a product is the Minkowski sum of the Newton polytopes of the
  factors. So consider the map $\PP^2 \to \PP^3$ defined by
  \[
    [ x_0^{} \mathbin{:} x_1^{} \mathbin{:} x_2^{}] \mapsto [ x_0^3
    \mathbin{:} x_1^2 x_2^{} \mathbin{:} x_1^{} x_2^2 \mathbin{:} x_0^{}
    x_1^{} x_2^{} ] \, .
  \]
  The image of this map is the zero-locus of the polynomial
  $q \coloneqq xyz - w^3$ where the homogeneous coordinate ring of $\PP^3$ is
  $\RR[x,y,z,w]$.  Hence, the Motzkin polynomial is the restriction of the
  function $f \coloneqq x^2 + y^2+ z^2 - 3w^2$ to the image variety.  As
  explained in the first section, this quadratic polynomial is not a sum of
  squares modulo the defining ideal of the surface because no nonzero
  polynomial of degree less than $3$ vanishes on the image variety and $f$ is
  not nonnegative on $\PP^3$. Therefore, the Motzkin polynomial cannot be a
  polynomial sum of squares.
\end{example}

The surprising answer to Question~\ref{q:general}, restated in the next
theorem, first appeared as Theorem~1.1 in \cite{BSV16}.  To formulate this
result, recall that a variety is irreducible if it is not the union of two
proper subvarieties. From the algebraic point of view, a variety is
irreducible if and only if the ideal of polynomials vanishing on it is prime.

\begin{theorem}[Blekherman, Smith, and Velasco]
  \label{thm:vmd}
  Let $X$ be an irreducible nondegenerate totally-real subvariety in
  $\PP^n$. Every nonnegative quadratic function on $X$ is a polynomial sum of
  squares, modulo the defining ideal of $X$, if and only if we have
  $\deg(X) = \codim(X) + 1$.
\end{theorem}

This theorem reveals a remarkable connection between a semi-algebraic property
(any nonnegative polynomial being a polynomial sum of squares) and the
fundamental geometric invariants of the complex variety.  For any subvariety
$X \subseteq \PP^n$, the codimension is a simple numerical measure of its
relative size; $\codim(X) \coloneqq n - \dim(X)$.  The degree is a second
numerical invariant depending on the embedding of the variety $X$ in $\PP^n$.
Geometrically, the degree, denoted by $\deg(X)$, is the number of points in
the intersection of the variety $X$ and a general linear subspace of dimension
$\codim(X)$.  For instance, if we intersect the twisted cubic
$C = \{[t_0^3 \mathbin{:} t_0^2 t_1^{} \mathbin{:} t_0^{} t_1^2 \mathbin{:}
t_1^3] \} \subset \PP^3$ with the generic hyperplane given by the linear
equation $a \, x_0 + b \, x_1 + c \, x_2 + d \, x_3 = 0$, then the
intersection points correspond to the three (typically distinct) roots of the
binary cubic
$a \, t_0^3 + b \, t_0^2 t_1^{} + c \, t_0^{} t_1^2 + d \, t_1^3$, so
$\deg(C) = 3$.  As the codimension of the curve $C$ in $\PP^3$ is $2$,
Theorem~\ref{thm:vmd} again implies that $P_C = \Sigma_C$; compare with
Example~\ref{e:twi}.

The irreducible nondegenerate subvarieties $X \subset \PP^n$ satisfying
$\deg(X) = \codim(X) + 1$ are called \emph{varieties of minimal degree}. As
the terminology suggests, these subvarieties do have the smallest possible
degree. The complete classification of varieties of minimal degree is one of
the outstanding achievements of the classical Italian school of algebraic
geometry.  Pasquale del Pezzo (1886) classified the surfaces of minimal degree
and Eugenio Bertini (1907) extended it to varieties of any dimension; see
\cite{MR927946} for an account.  The next subsection summarizes this result
and highlights links with sums of squares.  The remainder of this subsection
is devoted to two ideas: the impact of the minimal degree condition and the
pivotal role played by convex geometry, which has so far been hidden.

To get a feel for the degree condition, we temporarily assume that the
subvariety $X$ is a finite set of real points in $\PP^n$.  Although this
variety is reducible whenever there is more than one point, this
zero-dimensional case nevertheless develops some valuable intuition.  If $X$
is a variety of minimal degree, then it is a set of $n+1$ real points which
span the ambient space $\PP^n$.  Hence, we may choose coordinates on $\PP^n$
such that $X = \{[e_0], [e_1], \dotsc, [e_n] \}$ where $e_1, e_2, \dotsc, e_n$
are the standard basis vectors for $\RR^{n+1}$.

Why can we express every nonnegative quadratic form as a polynomial sum of
squares on this zero-dimensional variety?  The answer essentially comes from
interpolation.  The coordinate function $x_i$ vanishes at all points in $X$
except $[e_i]$ for all $0 \leqslant i \leqslant n$.  It follows that, for any
nonnegative quadratic form $q$ on $X$, the equality
\[
  q = \sum_{i = 0}^n \Bigl( \sqrt{q(e_i)} \, x_i \Bigr)^{\!2} 
\]
holds in the homogeneous coordinate ring of $X$, proving that $q$ is a
polynomial sum of squares.  If the subvariety $X$ has more than $n+1$ real
points (which implies that $\deg(X) > \codim(X) + 1$), then the interpolation
is no longer possible, because there exists a non-trivial linear relation
between the values of the linear forms at the points on $X$.  These linear
relations impose constraints on the possible values of quadratic sums of
squares on $X$ allowing us to separate $P_X$ and $\Sigma_X$.

For the surface $X = \nu_3(\PP^2) \subset \PP^9$, Hilbert~\cite{Hilbert1888}
already used these linear relations to prove that the cone $\Sigma_X$ of sums
of squares is properly contained in the cone $P_X$ of nonnegative polynomials.
A hyperplane section of the subvariety $X$ corresponds to a cubic curve in the
plane.  Cutting $X$ with two hyperplanes, one obtains the intersection of two
cubic curves in the plane.  By B\'ezout's Theorem, two general cubics
intersect in $9 = 3 \cdot 3$ points, so $\deg(X) = 9 = \codim(X) + 2$.  By
choosing the two cubics appropriately, we can arrange for the intersection
points to all be real.  These nine points of intersection lie in a $\PP^7$ and
are therefore linearly dependent. Moreover, the values of cubic forms on these
points are also linearly related. In classical algebraic geometry, this is
known as the Cayley--Bacharach theorem, which is usually stated as: any cubic
passing through any eight points of the intersection must also pass through
the ninth point.  Exploiting this linear relation, Hilbert showed that
sum-of-square cone $\Sigma_X$ is properly contained in the nonnegativity cone
$P_X$. Robinson later used Hilbert's technique to construct an explicit
nonnegative ternary form that is not a sum of squares; see \cite{Bruce}.

What serves as the bridge between complex algebraic geometry and
semi-algebraic nonnegativity results?  For the work under discussion, the
unifying answer is convex geometry.  Both $P_X$ and $\Sigma_X$ are more than
just subsets of the vector space $R_2$ of quadrics on the subvariety
$X \subseteq \PP^n$. They are closed convex cones. As convex cones, these
subsets are closed under taking linear combinations with nonnegative
coefficients.  Being closed means that these subsets are closed sets in the
natural Euclidean topology on the real vector space $R_2$.

Duality is an intrinsic feature of convex geometry.  Every closed convex cone
is equal to the intersection of the closed half-spaces that contain it.  The
linear inequalities defining these closed half-spaces form a convex cone in
the dual vector space, unimaginatively called the dual cone.  The most
accessible convex cones are polyhedral, defined by finitely many linear
inequalities.  Neither $P_X$ or $\Sigma_X$ is a polyhedral cone.  Their dual
cones can be quite complicated.  Fortuitously, the dual cone $\Sigma_X^*$
belongs to the next best class of convex cones.

For this part of our story, the cone of positive semidefinite matrices plays a
distinguished role. The positive semidefinite quadratic forms constitute the
closed convex cone $\mathcal{S}^{n+1}_+$ inside the real vector space
$\RR[x_0, x_1, \dotsc, x_n]_2^{}$ of quadratic polynomials.  Moreover, the
sum-of-squares cone $\Sigma_X$ is the image of the cone $\mathcal{S}^{n+1}_+$
under the canonical linear projection from $\RR[x_0, x_1, \dotsc, x_n]_2^{}$
to $R_2$.  It follows that the dual cone $\Sigma_X^*$ is
\emph{spectrahedral}---this cone can be represented as a linear matrix
inequality or equivalently as the intersection of $\mathcal{S}^{n+1}_+$ with a
linear affine subspace.  Understanding the minimal generators, also known as
extreme rays, of the spectrahedral cone $\Sigma_X^*$ underpins our
advancements.

Convex duality also produces fruitful connections between real algebraic
geometry and real analysis.  The dual cone $P_X^*$ consists of the linear
functionals that are nonnegative on $P_X$.  Given a linear functional that
evaluates nonnegatively on nonnegative functions, one hopes that it may be
represented as integration with respect to a measure supported on $X$.  This
forecast is often correct.  For example, the equality between nonnegative
polynomials and sums of squares in the univariate (or bivariate homogeneous
case) leads to a solution of the Hamburger, Hausdorff, and Stiltjes moment
problems; see Chapters~3 and~9 in \cite{momentproblem}.

\subsection*{The classification}

Varieties of minimal degree come in three flavors.  The most palatable are
quadratic hypersurfaces.  The real quadratic forms satisfying the hypothesis
of \Cref{thm:vmd} are necessarily indefinite.  Specializing to this case, the
theorem asserts that, when a quadratic form is nonnegative on the set of real
points at which an indefinite quadratic form vanishes, there exists a linear
combination of the two forms that is positive semidefinite.  This statement is
equivalent to the well-known $S$-Lemma (or $S$-procedure) in the optimization
community.  When $X = \PP^n$, this covers the quadratic forms in
Theorem~\ref{thm:hilb}.

The second flavor is an infinite family of projective varieties called
\emph{rational normal scrolls}.  Every member in this family is a smooth toric
variety, but there are infinitely many in each dimension.  The one-dimensional
family members are the rational normal curves arising as the Veronese
embeddings of $\PP^1$.  This family corresponds to the binary forms in
Theorem~\ref{thm:hilb}.

The third flavor consists of just one variety, namely the \emph{Veronese
  surface} in $\PP^5$.  This projective variety is defined by
$\{[x^2 \mathbin{:} xy \mathbin{:} xz \mathbin{:} y^2 \mathbin{:} yz
\mathbin{:} z^2]\in \PP^5 \mid [x \mathbin{:} y \mathbin{:} z]\in \PP^2\}$.
This exceptional variety corresponds to the sporadic case of the ternary
quartics in Theorem~\ref{thm:hilb}.

\begin{figure}[ht]
  \centering
  \begin{tikzpicture}[scale = 0.75]
    \filldraw[gray!50!white] (0,0) -- (4,0) -- (3,1) -- (0,1) -- (0,0);
    \draw[very thick] (0,0) -- (4,0) -- (3,1) -- (0,1) -- (0,0);
    \filldraw[black] (0,0) circle (3pt);
    \filldraw[black] (1,0) circle (3pt);
    \filldraw[black] (2,0) circle (3pt);
    \filldraw[black] (3,0) circle (3pt);
    \filldraw[black] (4,0) circle (3pt);
    \filldraw[black] (0,1) circle (3pt);
    \filldraw[black] (1,1) circle (3pt);
    \filldraw[black] (2,1) circle (3pt);
    \filldraw[black] (3,1) circle (3pt);
  \end{tikzpicture}
  \caption{The polygon of the rational normal scroll when $d = 4$ and
    $e = 3$}
  \label{fig:ratscroll}
\end{figure}
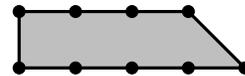

The structure of the rational normal scrolls warrants a closer look.  The
surfaces in this family are toric varieties corresponding to lattice polygons
of the form: $\operatorname{conv} \{ (0,0), (0,1), (d,0), (e,1) \}$; see
\Cref{fig:ratscroll}.
This polygon defines a map from the torus $(\CC^*)^2$ to $\PP^{d+e+1}$ given
by
$(s,t) \mapsto [1 \mathbin{:} t \mathbin{:} t^2 \mathbin{:} \dotsb \mathbin{:}
t^d \mathbin{:} s \mathbin{:} st \mathbin{:} \dotsb \mathbin{:} s t^e]$ where
monomials have exponent vectors equal to the lattice points of the polygon.
The closure of the image is the rational normal scroll corresponding to the
polygon.  It is a projective embedding of the $(d-e)$th Hirzebruch surface.
There are two rational normal curves, one of degree $d$ (for $s=0$) and one of
degree $e$ (for $s = \infty$) contained in this surface.  The scroll is swept
out by the lines joining the two points on the rational normal curves for the
same values of $t$.  Assuming that $d \geqslant e$, one may homogenize and
write as the map $\PP^1 \times \PP^1 \to \PP^{d+e+1}$ sending the pair
$( [s_0 \mathbin{:} s_1], [t_0 \mathbin{:} t_1] )$ to
\[
  [s_0^{} t_0^d \mathbin{:} s_0^{} t_0^{d-1} t_1^{} \mathbin{:} \dotsb
  \mathbin{:} s_0^{} t_1^d \mathbin{:} s_1^{} t_0^d \mathbin{:} s_1^{}
  t_0^{d-1} t_1^{} \mathbin{:} \dotsb \mathbin{:} s_1^{} t_0^{d-e} t_1^e ] \,
  .
\]
When $s_0 = 1$ and $t_0 = 1$, we recoup the first map.

Higher-dimensional rational normal scrolls are constructed in the same way. To
obtain a variety of dimension $k$, there are $k$ rational normal curves in
some $\PP^N$ (whose spans do not intersect) and the scroll is swept out by the
linear spaces spanned by the $k$ points on the curves for the same value of
$t$.

More generally, a lattice polytope determines a monomial map by interpreting
the lattice points lying in the polytope as exponent vectors. Toric geometry
provides a dictionary between the properties of the polytopes and properties
of the (closure of the) image of the corresponding monomial map.  A toric
surface is a variety of minimal degree if and only if the polygon does not
contain any lattice points in its interior.  Excluding the polygons of the
rational normal scrolls that have height one (like the one in
\Cref{fig:ratscroll}), there is only one more example: the triangle
$\conv\{(0,0),(2,0),(0,2)\}$ corresponding to the Veronese surface
$\nu_2(\PP^2) \subset \PP^5$.

Restricting a quadratic form to a rational normal scroll of dimension $2$
yields an element in two variables that has degree $2$ in $s$ and degree $2d$
in $t$ (say $d\geqslant e$). In the higher-dimensional cases, it is a form
that has degree $2$ in variables $s_1, s_2, \dotsc, s_{k-1}$ and some degree
$2d$ in $t$.  Homogenizing gives a form that is homogeneous of degree $2$ in
the $k$ variables $s_0, s_1, \dotsc, s_{k-1}$ and homogeneous of degree $2d$
in $t_0, t_1$.  In~\cite{CLR}, Man Duen Choi, TY Lam, and Bruce Reznick called
these element biforms.  A biform can also naturally be thought of as a
symmetric matrix with real polynomial entries that is pointwise positive
semidefinite on $\PP^1$ as in \cite{MR3968894}.

From a historical perspective, \Cref{thm:vmd} unites two fundamental results
that developed independently in the 1880s.  As we discuss the case of
reducible varieties, we will see that a similar story repeated about 100 years
later.

\subsection*{The reducible case}
\label{sec:reducible}

Thus far, we have emphasized sums of squares on irreducible varieties.
However, Theorem~\ref{thm:vmd} extends to all varieties including the
reducible ones. In this subsection, we discuss some of the ideas involved in
this generalization as well as some of its applications.

The varieties of minimal degree have the smallest possible degree among
irreducible varieties $X \subset \PP^n$ which span $\PP^n$.  In the reducible
case, this inequality no longer holds (think of two skew lines in
$\PP^3$). The right generalization of this geometric concept is not
immediately clear, but it turns out to be algebraic and involves syzygies.
Consider homogenous quadratic polynomials $f$ and $g$ in $3$ variables
defining two quadratic curves in $\PP^2$.  A \emph{syzygy} is a linear
relation between $f$ and $g$ where the coefficients are elements in the
polynomial ring $\RR[x_0,x_1,x_2]$.  If at least one of $f$ or $g$ is
irreducible, then every syzygy between them is generated by the obvious one:
$(g) f + (-f) g = 0$.  When $f$ and $g$ have a common linear factor. then a
linear syzygy exists.  For example, when $f = x_0 x_1$ and $g = x_0 x_2$, we
have $(x_2) f + (-x_1) g = 0$.

For a general set of polynomials, there can be syzygies among the syzygies.
This data is usually organized into a minimal free resolution.  One good way
to gauge the complexity of such a resolution is called the
Castelnuovo--Mumford regularity.  A linear space has Castelnuovo--Mumford
regularity $1$.  Varieties of minimal degree have the smallest
Castelnuovo--Mumford regularity among the varieties that are not linear
spaces, namely $2$.  For our purposes, $2$-regular varieties are the right
generalization for the concept of minimal degree.  The next result first
appeared as Theorem 9 in \cite{MR3633773}.

\begin{theorem}[Blekherman, Sinn, and Velasco]
  \label{thm:2regsos}
  Let $X$ be a nondegenerate totally-real projective variety in $\PP^n$. Every
  nonnegative quadratic form on $X$ is a sum of squares, modulo the defining
  ideal of $X$, if and only if $X$ is a $2$-regular variety.
\end{theorem}

What are all the $2$-regular varieties? In 2006, David Eisenbud, Mark Green,
Klaus Hulek, and Sorin Popescu gave a complete classification and a beautiful
geometric characterizations.  A subvariety $X \subseteq \PP^n$ is $2$-regular
if and only if, for any linear subspace $L$ such that $X \cap L$ is a finite
set, this finite set is linearly independent (linear independence has to be
defined carefully when the intersection is a non-reduced zero-dimensional
scheme).

Instead of digging into the details of this classification, we concentrate on
a special case, namely arrangements of coordinate subspaces.  We consider
projective varieties $X = \bigcup_{i=1}^k U_i$, where each $U_i$ is a linear
subspace in $\PP^{n-1}$ spanned by some set of standard coordinate vectors
$e_i$, where $1 \leqslant i \leqslant n$.  We also assume that corresponding
homogeneous ideals are generated in degree $2$.  For the combinatorially
inclined, these subspace arrangements correspond to flag complexes. The
definition ensures that each arrangement is determined by the coordinate lines
in $\PP^{n-1}$ that it contains.

\begin{figure}[ht]
  \centering
  \begin{tikzpicture}[scale = 0.6]
    \draw[very thick] (0,0) -- (-1.73205,1) -- (-1.73205,-1) -- (0,0);
    \draw[very thick] (0,0) -- (2,0);
    \filldraw[black] (-1.73205,1) circle (3pt) node[anchor=east] {$1$};
    \filldraw[black] (0,0) circle (3pt) node[anchor=south] {$3$};
    \filldraw[black] (-1.73205,-1) circle (3pt) node[anchor=east] {$2$};
    \filldraw[black] (2,0) circle (3pt) node[anchor=south] {$4$};
  \end{tikzpicture}
  \caption{The graph corresponding to the subspace arrangement in $\PP^3$
    defined by the ideal
    $\langle xw, yw \rangle = \langle w \rangle \cap \langle x, y \rangle$}
  \label{fig:graphMC}
\end{figure}
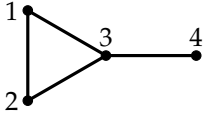

Alternatively, each such arrangement $X$ of coordinate subspaces is determined
by a graph.  The graph has $n$ vertices corresponding coordinates in
$\PP^{n-1}$.  There is an edge between vertex $i$ and vertex $j$ when
coordinate subspace spanned by the $i$th and $j$th variables is contained in
$X$. In other words, for each irreducible component $U_i \subset X$, we add
the clique on all the vertices $j$ satisfying $e_j \in U_i$. Hence, the ideal
defining the subvariety $X$ is generated by all the monomials $x_i x_j$ such
that $\{i,j\}$ is not an edge in our graph; see \Cref{fig:graphMC}.  To stress
that the subvariety $X$ comes from a graph $G$, we write $X(G)$.

The bijection between these subspace arrangements and graphs allows us to
translate combinatorial properties of the graph $G$ into geometric properties
of the subvariety $X(G)$. For instance, the dimension of $X(G)$ is one less
than the clique number $\omega(G)$ (the size of the largest clique in
$G$). For which graphs $G$ is $X(G)$ a $2$-regular variety? 

\begin{theorem}[Fr\"oberg] 
  \label{thm:frob}
  The variety $X(G)$ is $2$-regular if and only if the graph $G$ is chordal.
\end{theorem}

A graph $G$ is \emph{chordal} if it contains no induced cycles of length at
least four. In spite of their seemingly innocuous definition, chordal graphs
are very useful and arise naturally in combinatorics, numerical linear
algebra, and semidefinite programming; see \cite{AV}.
Theorem~\ref{thm:2regsos} shows that chordal graphs determine the only
varieties $X(G)$ on which nonnegative quadratic forms are sums of squares.

What does it mean for the quadratic form $q = x^\intercal \, A \, x$, where
$x \coloneqq \begin{bsmallmatrix*} x & y & z &
  w \end{bsmallmatrix*}^\intercal$, to be nonnegative when restricted to
$X(G)$?  The restriction remembers certain minors of the matrix $A$, but does
not know about other entries.  For the graph in \Cref{fig:graphMC}, the
restriction of the quatraic form $q$ to the $3$-dimensional linear subspace
defined by the ideal $\langle w \rangle$ corresponds to the upper left
$(3 \times 3)$-submatrix in $A$.  Similarly, the restriction of $q$ to the
$2$-dimensional linear subspace defined by the ideal $\langle x, y \rangle$
corresponds to the lower right $(2 \times 2)$-submatrix in $A$.  Hence, the
restriction has forgotten about the entries $a_{1,4}$ and $a_{2,4}$ in $A$. In
this sense, the restriction of the quadratic form $q$ to subspace arrangment
$X(G)$ is a partially specified real symmetric matrix. Nonnegativity of the
restriction means that the corresponding submatrices of $A$ are positive
semidefinite.

What does it mean for a nonnegative quadratic form on $X(G)$ to be a sum of
squares?  Our description of the ideal defining $X(G)$ implies that
restricting the quadratic form $x^\intercal \, A \, x$ to subvariety $X(G)$ is
the same as erasing the entries of $A$ corresponding to non-edges of $G$. To
write a quadratic form as a sum of squares modulo the ideal of $X$ is
tantamount to choosing the unspecified entries of the matrix so that the
entire symmetric matrix $A$ is positive semidefinite. This problem is known as
a positive semidefinite matrix completion problem.  It is well studied with
applications in combinatorics, discrete geometry, and statistics.

We can now reinterpret the equivalence between the chordality of $G$ and
$P_{X(G)} = \Sigma_{X(G)}$ in terms of matrix completion, reproving the
following theorem.

\begin{theorem}[Grone, Johnson, S\'a, and Wolkowicz]
  \label{thm:mc}
  Let $X(G)$ be the subspace arrangement associated to a graph $G$.  Every
  nonnegative quadratic form on $X(G)$ is a sum of squares, modulo the ideal
  of $X(G)$, if and only if the graph $G$ is chordal.
\end{theorem}

Remarkably, Theorem~\ref{thm:2regsos} proves that Theorems~\ref{thm:frob}
and~\ref{thm:mc} are in fact equivalent suggesting a link between matrix
completion problems and commutative algebra.

\section{Pythagoras numbers}
\label{sec:pythqp}

\noindent
In 1770, Joseph-Louis Lagrange proved his four-squares theorem: every
nonnegative integer can be represented as a sum of four integer squares.
Adrien-Marie Legendre extended this result in 1797 or 1798, proving that three
squares suffice unless the integer has the form $4^k(8m + 7)$ for some
nonnegative integers $k$ and $m$.  Pierre de Fermat had already determined in
1625 the exact conditions for a nonnegative integer to be expressible as a sum
of two squares.  Altogether this work gives a relatively complete description
of the sum-of-squares expressions for integers.  This achievement encourages
one to look for analogues for sums of squares in other rings.  This section
explores this line of inquiry for homogeneous coordinate rings of projective
real subvarieties.

Moving away from the relation between nonnegative elements and sums of
squares, we endeavor to stratify the various sum-of-squares decompositions of
an element according to their size.  The \emph{length} of an element $f$ in
the homogeneous coordinate ring of a projective subvariety is the smallest
number of summands needed to express $f$ as a sum of squares.

\begin{example}
  \label{e:length}
  The binary quartic $f \coloneqq (t_0^2)^{2} + (2 \, t_1^{2})^{2}$ has length
  two.  It has a second representation as the sum of two squares, namely
  $f = (t_0^2 - 2 \, t_1^2)^2 + (2 \, t_0^{} t_1^{})^2$.  Taking the average
  of these expressions, we obtain
  \[
    f
    = \Bigl( \!\tfrac{1}{\sqrt{2}}  t_0^2 \Bigr)^{\! 2}
    + \Bigl( \!\sqrt{2} \, t_1^2 \Bigr)^{\! 2}
    + \Bigl( \!\tfrac{1}{\sqrt{2}} t_0^{2} - \sqrt{2} \, t_1^2 \Bigr)^{\! 2}
    + \Bigl( \!\sqrt{2} \, t_0^{}  t_1^{} \!\Bigr)^{\! 2} \, ,
  \]
  and deduce that $f = (B t)^\intercal (B t)$ where
  \begin{align*}
    B
    &\coloneqq \frac{1}{\sqrt{2}}
      \begin{bmatrix*}[r]
        1 & 0 & 0 \\[-1pt]
        0 & 2 & 0 \\[-1pt]
        1 & -2 & 0 \\[-1pt]
        0 & 0 & 2 \\
      \end{bmatrix*}
    &&\text{and}
    & t &\coloneqq
      \begin{bmatrix*}
        t_0^2 \\[2pt] t_1^2 \\[1pt] t_0^{} t_1^{} \\[1pt]
      \end{bmatrix*} \, .     
  \end{align*}
  Under the orthogonal change of coordinates given by
  \[
    Q \coloneqq
    \frac{1}{6}
    \begin{bmatrix*}[r]
      3 \sqrt{2} & 0 & 3 \sqrt{2} & 0 \\[1pt]
      0 & 0 & 0 & 6 \\[1pt]
      \sqrt{6} & 2 \sqrt{6} & - \sqrt{6} & 0 \\[1pt]
      2 \sqrt{3} & -2 \sqrt{3} & -2 \sqrt{3} & 0 \\[1pt]
    \end{bmatrix*} \, ,
  \]
  this sum of four squares becomes a sum of three squares;
  $f = (Q B t)^\intercal (Q B t) = \smash{\bigl( t_0^2 \!-\! t_1^2 \bigr)}^2
  \!+\!  \smash{\bigl(\!\sqrt{2} \, t_0^{} t_1^{} \bigr)}^2 \!+\!
  \smash{\bigl( \!\sqrt{3} \, t_1^{2} \bigr)}^2$.
\end{example}

\begin{figure}[ht]
  \centering
    \includegraphics[width = \the\columnwidth]{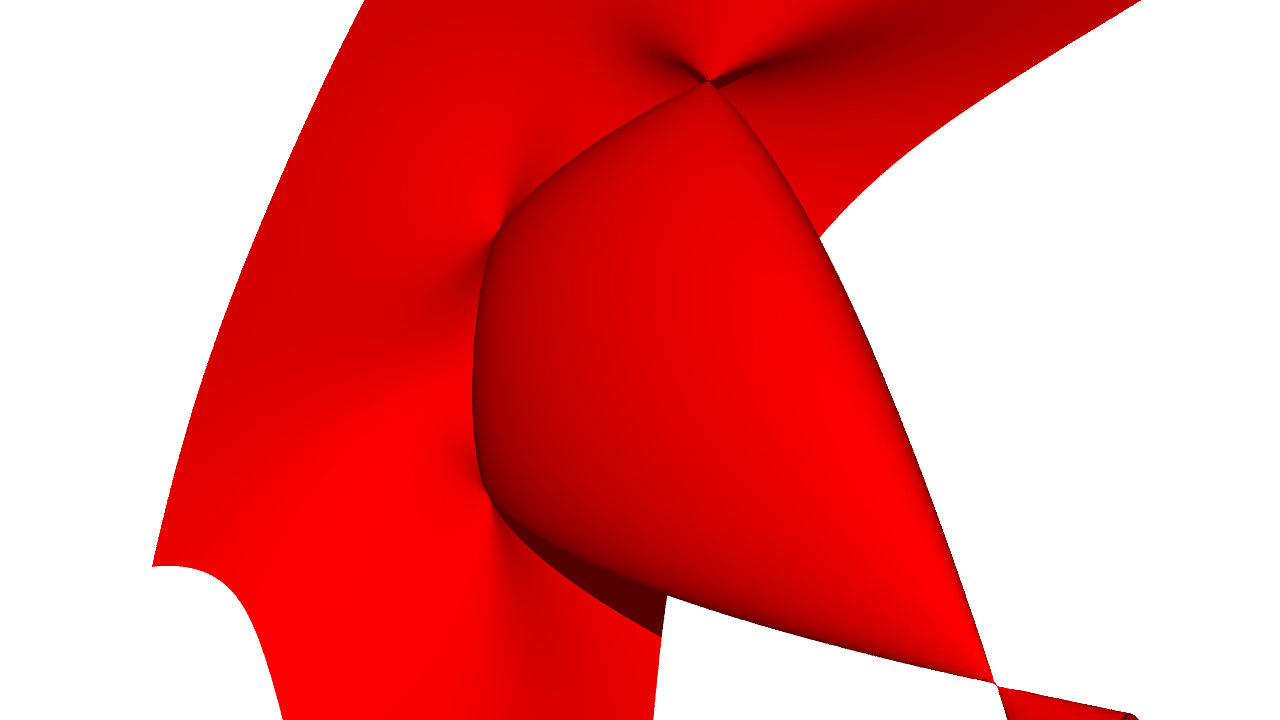}
  \caption{The Gram spectrahedron of a binary sextic}
  \label{f:spectrahedron}
\end{figure}

Sum-of-squares representations are parametrized by a convex semi-algebraic
set. The Gram spectrahedron of an element $f$ in the homogeneous coordinate
ring of a real subvariety $X \subseteq \PP^n$ consists of all positive
semidefinite quadratic forms in $\RR[x_0, x_1, \dotsc, x_n]$ that restrict to
$f$.  In Example~\ref{e:length}, the subvariety $X = \nu_2(\PP^1)$ is a conic
in $\PP^2$ and the Gram spectrahedron of the binary quartic is a line segment:
the convex hull of the two quadratic forms of rank $2$ whose restriction to
$\nu_2(\PP^1)$ is $f$.  The Gram spectrahedra of binary forms grow more
interesting as the degree increases.  For binary sextics,
Figure~\ref{f:spectrahedron} visualizes the Zariski closure of the boundary of
this $3$-dimensional Gram spectrahedron.  In this is case the boundary lies on
a Kummer surface and the spectrahedron is the convex region containing four
nodes; see Subsection~4.2 in \cite{Gramspecs}.

To see some repercussions of working modulo the ideal of a subvariety on
length, let $X$ be the hypersurface in $\PP^3$ defined by
$x^2 + y^2 + z^2 - w^2 = 0$. Consider the quadratic element in the homogeneous
coordinate ring of $X$ determined by $f = 2(x^2-yz) + w^2$.  Despite $f$
having a negative sign when evaluated at some points on $\PP^3$, the
corresponding element on $X$ is a sum-of-squares; homogenizing the first
equation in this article gives
\[
  2(x^2 \!-\! yz) + w^2 + ( x^2 \!+\! y^2 \!+\! z^2 \!-\! w^2)
  = 3x^2+ (y \!-\! z)^2 \, . 
\]
What is the length of $f$ on $X$?  The element $f$ vanishes at the real point
$p \coloneqq [0 \mathbin{:} 1 \mathbin{:} 1 \mathbin{:} \sqrt{2} ] \in
X$. Since a sum-of-squares representation of $f$ on $X$ evaluated at the point
$p$ is a sum of nonnegative real numbers adding to $0$, every summand also
vanishes at the point $p$.  It follows that every sum-of-squares
representation $f$ on $X$ is a sum of squares of linear forms vanishing at
$p$.  In more geometric terms, the representation for $f$ on $X$ factors
through a projection.  The projection away from the point $p$ is the rational
map $\pi_p \colon \PP^3 \dashrightarrow \PP^2$ whose components are given by a
basis of the linear forms vanishing at $p$ such as
\[
  \pi_p([x \mathbin{:} y \mathbin{:} z \mathbin{:} w])
  = \bigl[ x \mathbin{:} y - z \mathbin{:} \sqrt{2} y - w \bigr] \, .
\]
Writing $\RR[x_0,x_1,x_2]$ for the homogeneous coordinate ring of
$\pi_p(X) = \PP^2$, the element $f$ corresponds to $2x_0^2 + x_1^2$.  We
conclude that $f$ has length at least two on $X$.

To control the length of all elements, we introduce a new invariant.  For any
real subvariety $X \subseteq \PP^n$, the \emph{Pythagoras number} $\py(X)$ is
defined to be the smallest nonnegative integer $r$ such that any sum of
squares of linear elements in the homogeneous coordinate ring of $X$ can be
expressed as the sum of at most $r$ squares.  Although not apparent from the
definition, this semi-algebraic invariant reflects the underlying geometry of
the complex points in $X$.  For instance, we will see that having a small
Pythagoras number characterizes varieties of minimal degree. Additionally,
this numerical invariant is a crucial parameter in non-convex approaches to
semidefinite optimization, leading to valuable guarantees for sum-of-squares
methods on varieties; see \cite{MR1976484}.  Before delving into these
results, we revisit some examples from the previous section.

The simplest situation occurs when $X = \PP^n$. For any quadratic form $f$ on
$\PP^n$, there is a unique real symmetric matrix $A$ such that
$f = x^\intercal \, A \, x$ where $x \coloneqq
\begin{bsmallmatrix*}
  x_0 & x_1 & \dotsb & x_n \\
\end{bsmallmatrix*}^\intercal$.  The finite-dimensional spectral theorem
establishes that the matrix $A$ orthogonally diagonalizable.  This means that
there is an orthogonal matrix $Q$ such that $y \coloneqq Q \, x$ and
$f = \sum_{i=0}^n \lambda_i \, y_i^2$ where the real numbers
$\lambda_0, \lambda_1, \dotsc, \lambda_n$ are the eigenvalues of $A$.  The
quadratic form $f$ is nonnegative if and only if these eigenvalues are
nonnegative and have real square roots.  Thus, we may pull the coefficients
inside the squares and obtain an expression for $f$ as a sum of squares.  The
number of squares in this expression equals the number of nonzero eigenvalues,
so the length of the quadratic form $f$ is bounded above by the rank of
$A$. Conversely, if
$f = (v_1^\intercal x)^2 + (v_2^\intercal x)^2 + \dotsb + (v_k^\intercal x)^2$
is a sum of squares of linear forms whose coefficient vectors are
$v_1, v_2, \dotsc, v_k \in \RR^{n+1}$, then the corresponding matrix
$A = v_1^{} v_1^\intercal + v_2^{} v_2^\intercal + \dotsb + v_k^{}
v_k^\intercal$ is a sum of $k$ matrices, each having rank $1$.  Thus, the
matrix $A$ has rank at most $k$.  We surmise that the length of $f$ is equal
to the rank of $A$, which proves that $\py(\PP^{n}) = n+1$.  The same
reasoning implies that $\py(X) = n+1$ whenever the defining ideal of the real
subvariety $X \subseteq \PP^n$ contains no polynomials of degree less than
$3$.

Another situation is easy to understand.  Consider a nonnegative bivariate
form $f$ of degree $2j$, corresponding to a quadratic form on $\nu_j(\PP^1)$.
Dehomogenizing $f$ by setting $t_1 = 1$ yields a univariate polynomial
$q \in \RR[t_0]$.  Since the coefficients of $\hat{f}$ are real, any complex
roots come in conjugate pairs.  It follows that $f$ factors over the complex
numbers as $f = q_0 \, q_1 \, \overline{q_1}$ where $q_0$ is the product of
the linear forms corresponding to real roots of $\hat{f}$, $q_1$ is the
product of a linear form associated to each pair of complex roots of
$\hat{f}$, and $\overline{q_1}$ is the complex conjugate of $q_1$.  Since $f$
is nonnegative, every real root has even multiplicity, so $q_0$ is
automatically a square. From the identity
$q_1 \, \overline{q_1} = \operatorname{Re} (q_1)^2 + \operatorname{Im}
(q_1)^2$, we see that $f$ is a polynomial sum of two squares. In terms of
Pythagoras numbers, this shows $\py \bigl( \nu_j(\PP^1) \bigr) = 2$ for all
$j \geqslant 2$.

In our third situation, we focus on the appearance of Pythagoras numbers in
matrix completion problems.  As in the previous section, fix a pattern of
specified and unspecified off-diagonal entries in a matrix.  This data
determines a graph $G$ or square-free quadratic monomial ideal, whose
associated variety $X(G)$ is a coordinate subspace arrangement.  Rather than
looking for conditions that make a partially-specified matrix completable to a
positive semidefinite matrix, we now ask for low rank completions. Given a
matrix with missing entries that is known to have a positive semidefinite
completion, what is the lowest rank among its completions?  The smallest
nonnegative integer $r$, such that any completable partially-specified matrix
can be completed to a positive semidefinite matrix of rank $r$ is precisely
the Pythagoras number of variety $X(G)$.  In this context, the Pythagoras
number was called the Gram dimension by Monique Laurent and Antonios
Varvitsiotis in \cite{LV}. It has applications in distance realization
problems and rigidity theory.

Calculating Pythagoras numbers for the Veronese embeddings of projective space
is still an open problem.  For all positive integers $d$ and $n$, what is
$\py \bigl( \nu_d(\PP^n) \bigr)$?  To date, the strongest results are due to
Claus Scheiderer~\cite{MR3648509}.  He proves that
$\py \bigl( \nu_d(\PP^2) \bigr) \in \{d+1,d+2\}$. Assuming the validity of
conjectures by Anthony Iarrobino and Vassil Kanev on the Hilbert function of
ideals of finite sets of points \cite{MR1735271}, he also proves that the
asymptotic growth rate is $\py \bigl( \nu_d(\PP^n) \bigr) \in O(d^{n/2})$.
Our geometric perspective provides a systematic method for bounding Pythagoras
numbers.  The next lemma is a first step in this direction.

\begin{lemma}
  \label{lem: basic}
  Let $X \subseteq \PP^n$ be a totally-real subvariety. 
  \begin{compactenum}[\upshape (i)]
  \item For any totally-real subvariety $Y \subseteq \PP^n$ such that
    $X \subseteq Y$, we have the inequality $\py(X) \leqslant \py(Y)$.
  \item If $\pi_p \colon \PP^n \dashrightarrow \PP^{n-1}$ denotes the
    projection away from a real point $p \in X$, then we have
    $\py(X) \geqslant \py \bigl( \pi_p(X) \bigr)$. More precisely, the length
    of any quadratic element $f$ on $X$ that vanishes at the point $p$ is
    bounded below by the length of an element $g$ on $Y$ satisfying
    $f = \pi_p^*(g)$.
  \end{compactenum}
\end{lemma}

Part~(ii) showcases the importance of convexity.  By considering all elements
that vanish at a fixed real point $p \in X$, we identify a face $F$ of the
sum-of-squares cone $\Sigma_X$.  For any element $f$ lying on the face $F$,
every summand in any expression for $f$ as a sum of squares must also lie on
$F$.  It follows that the length of an element in $F$ is determined by the
face alone.  Moreover, the face $F$ of $\Sigma_X$ specified by evaluation at
the point $p \in X$ is isomorphic to the cone $\Sigma_Z$ where
$Z \coloneqq \pi_p(X)$.  Hence, the Pythagoras number of $X$ is bounded below
by the Pythagoras number of $\pi_p(X)$.

To estimate the Pythagoras number of a totally-real subvariety
$X \subseteq \PP^n$, Lemma~\ref{lem: basic} suggests that we recognize a
simpler subvariety containing $X$ or successively project away from points to
obtain a simpler variety.  For either approach to work, we need a sufficiently
large class of varieties having a known Pythagoras number.  Fortunately,
Corollary~32 in \cite{MR3633773} and Theorem~2.1 in \cite{BPSV} demonstrate
that varieties of minimal degree and their reducible generalizations, called
$2$-regular varieties, form such a class of varieties.

\begin{theorem}[Blekherman, Sinn, and Velasco]
  Let $X$ be a nondegenerate totally-real $2$-regular subvariety in $\PP^n$.
  We have $\py(X) = \dim (X) + 1$.  When $X$ is irreducible, we have
  $\py(X) = \dim(X) + 1$ if and only if $X$ is a variety of minimal degree.
\end{theorem}

\noindent
We may bound the Pythagoras number of a totally-real subvariety by embedding
it in a $2$-regular variety.

This geometric approach applies to matrix completion problems.  Recall that a
$2$-regular quadratic square-free monomial ideal corresponds to a chordal
graph $H$ and the dimension of the associated variety $X(H)$ is equal to the
size of the largest clique of $H$ minus one.  By considering all chordal
graphs that contain a given graph $G$, we have
$\py \bigl( X(G) \bigr) \leqslant \operatorname{tw}(G) + 1$, where
$\operatorname{tw}(G)$ is the tree-width of graph $G$; as in \cite{LV}, the
tree-width is defined to be one less than the minimum clique number among all
chordal graphs containing $G$.

\subsection*{Quadratic persistence}

As we have already observed, the equation $\py(Z) = m + 1$ holds whenever the
defining ideal of the real subvariety $Z \subseteq \PP^m$ contains no
polynomials of degree less than $3$. Using this observation and
Lemma~\ref{lem: basic}, we may bound the Pythagoras number of a totally-real
subvariety $X \subseteq \PP^n$ by successively projecting it away from $k$
real points so that no quadric polynomial vanishes on the image variety
$Z \subseteq \PP^{n-k}$.  If it is possible to eliminate all quadrics using
$k$ projections, then one would obtain $\py(X) \geqslant \py(Z) = n-k+1$. The
cardinality of a smallest set of points on $X$ which, via projection,
eliminate all quadrics is called the \emph{quadratic persistence} of $X$ and
denoted $\qp(X)$.  When $X$ is irreducible, choosing points at random will
always suffice.  This numerical invariant is purely algebraic because it
measures vanishing of quadratic forms (including the complex quadratic forms
in the ideal of $X$) on sets of points on the variety.

To become more acquainted with this process, Table~\ref{table} catalogues the
dimension of the vector space of quadrics polynomials vanishing on various
projections of $\nu_d(\PP^n)$ for small values of $d$ and $n$. The $k$th row
displays the dimension of the vector space of quadrics vanishing in the
projection away from $k$ generic points.
\begin{table}[ht]
  \caption{Dimensions of the space quadrics}
  \label{table}
  \begin{tabular}{ccccc}
    $\#$ Points & $\nu_2(\PP^2) \subset \PP^5$
    & $\nu_3(\PP^2) \subset \PP^9$
    & $\nu_3(\PP^3) \subset \PP^{19}$ \\ \hline
    0 & 6 & 27 & 126 \\
    1 & 3 & 20 & 110 \\
    2 & 1 & 14 & 95 \\
    3 & 0 & 9 & 81 \\
    4 & 0 & 5 & 68 \\
    5 & 0 & 2 & 56 \\
    6 & 0 & 0 & 45 \\
  \end{tabular}
\end{table}
This table shows that quadratic persistences of the varieties in second and
third columns are equal to $3$ and $6$ respectively.  It follows that
$\py \bigl( \nu_2(\PP^2) \bigr) \geqslant 3$ and
$\py \bigl( \nu_3(\PP^2) \bigr) \geqslant 4$.  Both of these inequalities are,
in fact, sharp.  The last column only shows that
$\qp \bigl( \nu_3(\PP^3) \bigr) \geqslant 7$.

How many quadrics do we lose in each step? It is straightforward to show that,
in each projection step, the number of quadrics that disappear from the ideal
is at most the codimension of the variety.  Surprisingly, this inequality is
an equality for each entry in Table~\ref{table}.  If this pattern were to
continue, then we would deduce that $\qp \bigl( \nu_3(\PP^3) \bigr) = 12$. The
Iarrobino--Kanev conjectures mentioned earlier imply that the number of
quadrics lost at every step is maximal for all values of $n$ and $d$.  Thus,
they imply $\qp \bigl( \nu_3(\PP^3) \bigr) = 12$.

More generally, for an irreducible subvariety $X \subseteq \PP^n$, the
projection away from a set of $\codim(X)$ generic points maps $X$ surjectively
onto $\PP^{\dim(X)}$.  Hence, there are no quadric polynomials vanishing on
the image.  It follows that $\qp(X) \leqslant \codim(X)$ and
$\py(X) \geqslant \dim(X) + 1$.

Quadratic persistence depends only on the complex algebraic geometry of a
variety $X$. It is an algebraic invariant rather than a semi-algebraic one, so
tools from commutative and homological algebra apply.  Describing this in
detail would take us too far afield.  Nevertheless, we want to state two
consequences of these tools, namely the characterizations
(modulo a few technical assumptions) of all irreducible varieties with minimal
and next-to-minimal Pythagoras numbers.  The following theorems appears as
Theorems 1.4 and 1.5 in \cite{BSV20}.

\begin{theorem}[Blekherman, Smith, Sinn, and Velasco]
  Let $X \subseteq \PP^n$ be a nondegenerate totally-real irreducible
  subvariety. The following three conditions are equivalent:
  \begin{compactenum}[\upshape (1)]
  \item $\py(X) = \dim(X) + 1$,
  \item $\deg(X) = \codim(X) + 1$,
  \item $\qp(X) = \codim(X)$.
  \end{compactenum}
\end{theorem}

\begin{theorem}[Blekherman, Smith, Sinn, and Velasco]
  Let $X \subseteq \PP^n$ be a nondegenerate totally-real irreducible
  subvariety. If $X$ is arithmetically Cohen-Macaulay, then the following
  three conditions are equivalent:
  \begin{compactenum}[\upshape (1)]
  \item $\py(X) = \dim(X) + 2$,
  \item $\deg(X) = \codim(X) + 2$ or $X$ is a subvariety having codimension
    $1$ in a variety of minimal degree,
  \item $\qp(X) = \codim(X) - 1$.
  \end{compactenum}
\end{theorem}

In all these cases and in all cases that we understand, quadratic persistence
detects the Pythagoras number: $\qp(X) + \py(X) = n+1$.  We have no reason to
expect this behavior, so it would be interesting to find an example where this
equation does not hold.

\section{Sum-of-squares multipliers}
\label{sec:multipliers}

\noindent
Knowing that there exist nonnegative polynomials that cannot be represented as
a polynomial sum of squares, we want other effective ways of recognizing
nonnegativity. The 17th problem on Hilbert's celebrated list proposes a
candidate.  Hilbert asks whether every polynomial that is nonnegative, when
regarded as a function from $\RR^n$ to $\RR$, equals a sum of squares of
rational functions.

Emil Artin~\cite{Artin} solved this problem in 1927.  Before stating this
result, observe that the intermediate value theorem implies that every real
polynomial having odd degree is negative at some point in $\RR^n$.  Moreover,
the homogenization of a polynomial takes only nonnegative values if and only
if the same is true for the dehomogenized polynomial.  Hence, we may restrict
our attention to homogeneous polynomials of even degree.

Given a nonnegative polynomial $f$, Artin demonstrates that, for some positive
integer $r$, there exist polynomials $g_1, g_2, \dotsc, g_{r}$,
$h_1, h_2,\dotsc, h_{r}$ such that
\[
  f = \left( \frac{h_1}{g_1} \right)^{\!\! 2} + \left(
    \frac{h_2}{g_2} \right)^{\!\! 2} + \dotsb + \left(
    \frac{h_{r}}{g_{r}} \right)^{\!\! 2} \, .
\]
For instance, the Motzkin polynomial is equal to a sum of four squares of
rational functions:
\begin{align*}
  &\relphantom{=}
    x_{0}^6 + x_{1}^4 x_{2}^2 + x_{1}^2 x_{2}^4 - 3 x_{0}^2 x_{1}^2 x_{2}^2\\
  &= \left(\! \frac{x_{0}x_{1}x_{2}(x_{1}^2 + x_{2}^2 - 2x_{0}^2)}%
    {x_{1}^2 + x_{2}^2} \!\right)^{\!\! 2}
    + \left(\! \frac{x_{1}^2x_{2}(x_{1}^2 + x_{2}^2 - 2x_{0}^2)}%
    {x_{1}^2 + x_{2}^2} \!\right)^{\!\! 2} \\
  &\relphantom{=W}
    + \left(\! \frac{x_{1}x_{2}^2(x_{1}^2 + x_{2}^2 - 2x_{0}^2)}%
    {x_{1}^2 + x_{2}^2} \!\right)^{\!\! 2}
    + \left(\! \frac{x_{0}^3(x_{1}^2 - x_{2}^2)}%
    {x_{1}^2 + x_{2}^2} \!\right)^{\!\! 2} \, . 
\end{align*}

Artin's original solution comes from studying sums of squares in arbitrary
fields.  The key insight uses total orderings that are compatible with field
operations to characterize the subset of elements in a field that can be
expressed as a sum of squares.  Although profoundly influential in the theory
of real-closed fields, this nonconstructive proof does not bound the number of
squares or the degrees of the polynomials appearing in the rational functions.
Albrecht Pfister~\cite{Pfister} subsequently proved that, for a homogeneous
polynomial in $n$ variables, $2^n$ squares always suffice.

In our search for constructive methods of identifying nonnegative polynomials,
we prefer a slightly different formulation.  We abandon the rational functions
by finding a common denominator.  Given a homogenous polynomial $f$ that is
nonnegative on $\RR^n$, it follows that, for some positive integers $s$ and
$r$, there are homogeneous polynomials
$g_1, g_2, \dotsc, g_{s}, h_1, h_2, \dotsc, h_{r}$ such that
\[
  (g _1^2 + g_2^2 + \dotsb + g_{s}^2) f = h_1^2 + h_2^2 + \dotsb + h_{r}^2
  \, .
\]
Allowing the multiplier $g \coloneqq g_1^2 + g_2^2 + \dotsb + g_{s}^2$ to be
any polynomial sum of squares rather than just the square of a single
homogeneous polynomial enlarges the pool of potential certificates.  At the
expense of doubling the degree of the multiplier and, at worst, replacing $r$
by $rs$, we recover a representation of $f$ as a sum of squares of rational
functions from $g^2 f = g (h_1^2 + h_2^2 + \dotsb + h_{r}^2)$.  To effectively
certify nonnegativity, one seeks to bound the degree on a multiplier $g$.

An explicit bound on multipliers is a recent discovery.  Henri Lombardi,
Daniel Perrucci and Marie-Fran\c{c}oise Roy~\cite{LPR} prove that, for any
nonnegative polynomial $f$ of degree $d$ in $n$ variables, there exists a
sum-of-squares multiplier $g$ of degree less than
\[
  \text{\large 2}^{\text{\normalsize 2}^{\text{\small 2}^%
      {\text{\footnotesize $d$\hspace{0.6pt}}^%
        {\scriptstyle 4^{\scriptscriptstyle n}}}}} \, .
\]
This tower of five exponents arises from bounding the complexity of a
quantifier elimination problem.  Roughly speaking, a cylindrical algebraic
decomposition produces two levels, a constructive version of the fundamental
theorem of algebra gives one level, and constructions based on the
intermediate value theorem are responsible for the other two levels in the
tower.

In contrast, very little is known about worst-case lower bounds on the degree
of a multiplier.  Grigoriy Blekherman, Jo\~ao Gouveia, and James
Pfeiffer~\cite{BGP} establish that, for all positive integers $n$, there is a
nonnegative quartic polynomial $f$ in $n$ variables such that any
sum-of-squares multiplier $g$ must have degree at least $n$.  Thus, we see
that there is an embarrassingly large gap between the current upper and lower
bounds on the degree of a multiplier.  In an attempt to redress this issue, we
aim to provide much tighter degree bounds in some situations.  We again appeal
to complex projective geometry.

Unlike in the previous sections, we do not limit ourselves to quadratic
functions on a real subvariety $X \subseteq \PP^n$. As before, let $R$ be the
homogeneous coordinate ring of $X$. For any nonnegative integer $j$, let
$P_{X,2j}$ denote the set of homogeneous elements in $R$ of degree $2j$ that
are nonnegative on $X$.  Similarly, let $\Sigma_{X,2j}$ denote the set of
homogeneous elements in $R$ of degree $2j$ that are sums of squares.  Compared
with our earlier notation, $P_{X,2j}$ and $\Sigma_{X,2j}$ are the same as
$P_{\nu_j(X)}$ and $\Sigma_{\nu_j(X)}$.  Hence, both $P_{X,2j}$ and
$\Sigma_{X,2j}$ are full-dimensional convex cones in the real vector space
$R_{2j}$ that contain no line and are closed in the Euclidean
topology. Considering higher-degree elements in $R$ seems more natural for
multipliers.

The existence of a nonnegative multiplier $g$ such that the product $gf$ is a
sum of squares manifestly confirms that the element $f$ is nonnegative when
evaluated at any point for which $g$ does not vanish.  The consequences for
other points are not immediately as clear.  When the complement of this
vanishing set is dense in the Euclidean topology, it follows that the element
$f$ is nonnegative at every point.  This rationale always applies to
polynomials on $\RR^n$, but becomes more complicated when the algebraic
variety $X$ is reducible or singular.

For real projective curves, we produce sharp degree bounds on sum-of-squares
multipliers interms of three fundamental numerical invariants.  As in
Section~\ref{sec:psdsos}, the degree of the curve $X \subset \PP^n$ is the
number points in the intersection of the curve and a general
$(n-1)$-dimensional linear subspace.  The (arithmetic) genus of a projective
curve $X$ is denoted by $\genus(X)$.  When $X$ is nonsingular over $\CC$, this
complex curve may be viewed as a Riemann surface and this numerical invariant
coincides with its topological genus.  Our third invariant $\hilbreg(X)$ is
the smallest integer $j$ such that the Hilbert function and the Hilbert
polynomial of $X$ agree when evaluated at any integer greater than or equal to
$j$.  Using these invariants, Theorem~1.1 in \cite{BSV20} gives the following
bounds.

\begin{theorem}[Blekherman, Smith, and Velasco]
  \label{t:curve}
  For any nondegenerate totally-real projective curve $X \subset \PP^n$, any
  nonnegative element $f \in P_{X,2j}$, and all nonnegative integers $k$
  satisfying $k \geqslant \max \{ 2 \genus(X)/\deg(X), \hilbreg(X) \}$, there
  exists a nonzero sum of squares $g \in \Sigma_{X,2k}$ such that the product
  $g f \in \Sigma_{2(j+k)}$ is also a sum of squares.

  Conversely, for all integers $n$ and $j$ greater than $1$, there exists a
  totally-real smooth curve $X \subset \PP^n$ and a nonnegative element
  $f \in P_{X,2j}$ such that, for all nonnegative integers $k$ satisfying
  $k < \max\{ 2 \genus(X)/ \deg(X), \hilbreg(X) \}$ and all nonzero sum of
  squares $g \in \Sigma_{X, 2k}$, the product $g f \not\in \Sigma_{X,2(j+k)}$
  is not a sum of squares.
\end{theorem}

The uniform degree bound on the multiplier $g$ in the first part of
Theorem~\ref{t:curve} is determined by just the complex geometry of the curve
$X$.  It is, notably, independent of both the degree of the nonnegative
element $f$ and the topology of the real points in $X$.

\begin{example}
  Consider a totally-real complete intersection curve $X \subset \PP^n$ cut
  out by homogeneous polynomials of degrees $d_1, d_2, \dotsc, d_{n-1}$.  The
  numerical invariants of interest are $\deg(X) = d_1 d_2 \dotsb d_{n-1}$,
  \[
    \genus(X) = \tfrac{1}{2}d_1 d_2 \dotsb d_{n-1}
    (d_1 \!+\! d_2 \!+\! \dotsb \!+\! d_{n-1} \!-\! n \!-\! 1) + 1
  \]
  and
  $\hilbreg(X) = d_1 \!+\! d_2 \!+\! \dotsb \!+\! d_{n-1} \!-\! n \!-\! 1$.
  Theorem~\ref{t:curve} shows that, for all
  $k \geqslant d_1 \!+\! d_2 \!+\! \dotsb \!+\! d_{n-1} \!-\! n$ and all
  nonnegative elements $f \in P_{X,2j}$, there is a nonzero sum of squares
  $g \in \Sigma_{X,2k}$ such that the product $gf \in \Sigma_{X,2(j+k)}$ is a
  sum of squares.
\end{example}

\begin{example}
  Theorem~\ref{t:curve} is not sharp on every curve.  Let $X$ be the image of
  the map $\PP^1 \to \PP^2$ defined by
  $[t_{0}^{} \mathbin{:} t_{1}^{}] \mapsto [ t_{0}^2
  t_{1}^{}(t_{0}^{}-t_{1}^{}) \mathbin{:} t_{0}^{} t_{1}^2(t_{0}^{} -
  t_{1}^{}) \mathbin{:} t_{0}^4+t_{1}^4]$.  This planar curve has degree $4$,
  arithmetic genus $3$, and $\hilbreg(X) = 2$, so Theorem~\ref{t:curve} would
  require $k \geqslant 2$.  However, one can prove that, for all nonnegative
  elements $f \in P_{X,2}$, there exists a nonzero sum of squares
  $g \in \Sigma_{X,2}$ such that the product $g f \in \Sigma_{X,4}$ is a sum
  of squares.  For more information about this analysis, see Example~5.3 in
  \cite{BSV20}.
\end{example}

The proofs for the two parts of Theorem~\ref{t:curve} are disjoint.  The upper
bound on the minimum degree of a multiplier is derived from a new Bertini
theorem in convex algebraic geometry and the lower bound is obtained by
deforming rational Harnack curves on toric surfaces.

For the first part, we reinterpret the non-existence of a sum-of-squares
multiplier $g \in \Sigma_{X, 2k}$ as asserting that the convex cones
$\Sigma_{X,2j+2k}$ and $f \cdot \Sigma_{X,2k}$ intersect only at zero.  If a
real subvariety $X \subseteq \PP^n$ has a linear functional separating these
cones, then we show that a sufficiently general hypersurface section of $X$
also does.  The phrase ``sufficiently general'' means belonging to a nonempty
open subset in the Euclidean topology of the relevant parameter space.
Unexpectedly, this convex variant of the Bertini Theorem relies on a
characterization of spectrahedra that have many facets in a neighbourhood of
every point.  Recognizing this dependency is the main insight.  By repeated
applications of our Bertini theorem, we reduce to the case of points.

For the second part, we establish that having a nonnegative element vanish at
a relatively large number of isolated real singularities prevents it from
having a low-degree sum-of-square multiplier.  The hypotheses needed to
realize this basic premise are formidable.  Nonetheless, this transforms the
problem into finding enough curves that satisfy the conditions and maximizing
the number of isolated real singularities.  Harnack curves having the maximal
number of connected real components in a prescribed topological arrangement,
are natural contenders.  We confirm that rational singular Harnack curves on
toric surfaces fulfill the technical requirements.  By perturbing both the
curves $X$ and the nonnegative element $f \in R_{2j}$, we exhibit smooth
curves and nonnegative elements without low-degree sum-of-squares multipliers.
Miraculously, for totally-real projective curves, the upper bounds and lower
bounds on the degree coincide.

This approach also works for surfaces of minimal degree.  Recall that a
subvariety $X \subset \PP^n$ has minimal degree if it is nondegenerate and
$\deg(X) = 1 + \codim(X)$.  Even when $X$ is a variety of minimal degree, its
$j$th Veronese embedding $\nu_j(X)$ rarely is.  Hence, multipliers are need to
certify nonnegativity.  Theorem~1.2 in \cite{BSV20} provides the following
bounds.

\begin{theorem}[Blekherman, Smith, and Velasco]
  \label{t:surface}
  Let $X$ be a totally-real surface of minimal degree in $\PP^n$.  For any
  nonnegative element $f \in P_{X,2j}$, there exists a nonzero sum of squares
  $g \in \Sigma_{X, j^2 - j}$ such that the product $gf \in \Sigma_{X,j^2+j}$
  is a sum of squares.

  Conversely, for any integer $j$ greater than $4$, there exists a nonnegative
  element $f \in P_{X, 2j}$ such that, for any positive integer $k$ satisfying
  $k < j-2$ and any nonzero sum of squares $g \in \Sigma_{X,2k}$, the product
  $g f \not\in \Sigma_{X,2(j+k)}$ is not a sum of squares.
\end{theorem}

Unlike for curves, Theorem~\ref{t:surface} demonstrates that the minimum
degree of a multiplier $g$ generally depends on the degree of the nonnegative
element $f$.  Furthermore, our techniques do not typically yield sharp bounds.

Even in the special case $X = \PP^2$, our geometric outlook leads to something
new.  For example, we re-prove and prove the following two results for ternary
octics:
\begin{compactitem}[$\bullet$]
\item For all nonnegative elements $f \in P_{\PP^2, 8}$, there exists a
  nonzero sum of squares $g \in \Sigma_{\PP^2, 4}$ such that the product
  $g f \in \Sigma_{\PP^2, 12}$ is also a sum of squares.
\item There exists a nonnegative element $f \in P_{\PP^2, 8}$ such that, for
  all nonzero sum of squares $g \in \Sigma_{\PP^2,2}$, the product
  $g f \not\in \Sigma_{\PP^2,10}$ is not a sum of squares.
\end{compactitem}
Together these observations give the first tight bounds on the degrees of
sum-of-squares multipliers for polynomials since Hilbert's 1893
paper~\cite{Hilbert1893} where he establishes sharp bounds for ternary
sextics.

\section*{Continuing the story}

\noindent
This is not the end of the story.  We have no doubt that many fascinating new
chapters about the connections between complex algebraic geometry and
nonnegativity have yet to be written.  In particular, we would like to see
solutions to the following open problems.
\begin{compactitem}[$\bullet$]
\item[\emph{Quantifying the discrepancy}:] How should one measure the
  difference between the nonnegative cone $P_X$ and the sum-of-squares cone
  $\Sigma_X$ when the subvariety $X \subset \PP^n$ is not of minimal degree?
  For some estimates on their relative volumes, see Chapter~4 in
  \cite{SIAMbook}.
\item[\emph{Pythagoras numbers}:] Can one clarify the relation between
  quadratic persistence and the Pythagoras number? For instance, is the
  classification of the real subvarieties $X \subset \PP^n$ satisfying
  $\qp(X) + \py(X) = n+1$ in sight?
\item[\emph{Sharp degree bounds}:] Which geometric invariants of the
  associated complex projective variety determine the smallest degree of a
  sum-of-squares multiplier on a surface or higher-dimensional variety?
\end{compactitem}
We recommend the recent books \cite{SIAMbook} and \cite{AMSshort} to readers
interested in further exploring these questions.  The authors thank the SIAM
Activity Group on Algebraic Geometry whose welcoming and interdisciplinary
community has played an important role in facilitating our
collaboration. Please join us in discovering more intriguing connections!

\subsection*{Acknowledgments}

\noindent
Grigoriy Blekherman was partially supported by the National Science Foundation
(NSF) grant DMS-1901950, Gregory G. Smith was partially supported by the
Natural Sciences and Engineering Research Council of Canada (NSERC), and
Mauricio Velasco was partially supported by proyecto INV-2018-50-1392 from
Facultad de Ciencias, Universidad de los Andes.


\raggedright

\end{document}